\documentclass{gtart_h}  

\def\ifplaintex{\expandafter\ifx\csname documentclass\endcsname\relax}

\ifplaintex 
\else
\fi

\expandafter\ifx\csname beginpicture\endcsname\relax
\expandafter\ifx\csname documentclass\endcsname\relax
\input pictex \else
\input prepictex \input pictex \input postpictex \fi\fi

\def\gt{{\mathsurround=0pt\it $\cal G\mskip-2mu$eometry \&\ 
$\cal T\!\!$opology}}        

\def\gtp{{\mathsurround=0pt\it $\cal G\mskip-2mu$eometry \&\ 
$\cal T\!\!$opology $\cal P\!$ublications}}  

\def\lognumber#1{\def\thelognumber{#1}}
\def\volumenumber#1{\def\thevolumenumber{#1}}
\def\papernumber#1{\def\thepapernumber{#1}}
\def\volumeyear#1{\def\thevolumeyear{#1}}

\def\pagenumbers#1#2{\def\startpage{#1}\def\finishpage{#2}}
\def\published#1{\def\publishdate{#1}}
\def\proposed#1{\def\theproposer{#1}}
\def\seconded#1{\def\theseconders{#1}}
\def\received#1{\def\receiveddate{#1}}

\def\accepted#1{\def\accepteddate{#1}}

\def\asciiaddress#1{\def\theasciiaddress{#1}}

\long\def\asciiabstract#1{\long\def\theasciiabstract{#1}}
\def\asciikeywords#1{\def\theasciikeywords{#1}}

\let\\\par\let\thelognumber\relax
\let\thevolumenumber\relax\let\thepapernumber\relax
\let\thevolumeyear\relax\let\thesamplenumber\relax\let\startpage\relax
\let\finishpage\relax\let\publishdate\relax\let\receiveddate\relax
\let\reviseddate\relax\let\accepteddate\relax\let\theasciititle\relax
\let\theasciiauthors\relax\let\theasciiaddress\relax
\let\theasciiabstract\relax\let\theasciikeywords\relax
\let\theasciiemail\relax\let\theshortauthors\relax\let\theshorttitle\relax

\long\def\maketitlep{   

\count0=\startpage

\gt\hfill      
\beginpicture
\setcoordinatesystem units <0.33truein, 0.33truein> point at 2.2 0.9
\setplotsymbol ({$\cal G$})
\plotsymbolspacing=9truept
\circulararc 315 degrees from 0 1 center at 0 0
\setplotsymbol ({$\cal T$})
\circulararc 315 degrees from 1 -1 center at 1 0
\endpicture
\break
{\small\ifx\thesamplenumber\relax 
Volume \else Sample
\fi\thevolumenumber\ (\thevolumeyear)
\startpage--\finishpage\nl
Published: \publishdate}
\vglue 0.5truein plus 0.4fil minus 0.1truein

{\parskip=0pt\leftskip 0pt plus 1fil\def\\{\par\smallskip}{\ifplaintex\large
\else\Large\fi\bf\thetitle}\par\medskip}   

\vglue 0pt plus 0.1fil 

{\parskip=0pt\leftskip 0pt plus 1fil\def\\{\par}{\sc\theauthors}
\par\medskip}

\vglue 0pt plus 0.1fil 

{\small\parskip=0pt\let\newline\\
{\leftskip 0pt plus 1fil\def\\{\par}{\sl\theaddress}\par}
\expandafter\ifx\theemail\relax    
\relax\else\vglue 5pt plus 0.02fil minus 2pt\def\\{\stdspace{\rm 
and}\stdspace} 
\cl{Email:\stdspace\tt\theemail}\fi
\ifx\theurl\relax                  
\relax\else\vglue 5pt plus 0.02fil minus 2pt\def\\{\stdspace{\rm 
and}\stdspace}
\cl{URL:\stdspace\tt\theurl}\fi\par}

\vglue 7pt plus 0.3fil minus 3pt

{\bf Abstract}
\vglue 5pt plus 0.1fil minus 2pt

\theabstract

\vglue 7pt plus 0.3fil minus 3pt

{\bf AMS Classification numbers}\quad Primary:\quad \theprimaryclass

Secondary:\quad \thesecondaryclass

\vglue 5pt plus 0.3fil minus 2pt

{\bf Keywords:}\quad \thekeywords

\vglue 10pt plus 0.5fil minus 5pt

{\small  Proposed: \theproposer\hfill Received: \receiveddate\nl
Seconded: \theseconders\hfill 
\ifx\reviseddate\relax                         
Accepted: \accepteddate                        
\else
Revised: \reviseddate                          
\fi}
\eject
}       

\let\maketitlepage\maketitlep
\let\maketitle\maketitlepage

\font\phead=cmsl9 scaled 950
\font=cmsl9 scaled 1050
\font\pnum=cmbx10 scaled 913
\font\lnum=cmbx10 
\font\pfoot=cmsl9 scaled 950
\font=cmsl9 scaled 1050
\ifplaintex
\headline{\vbox to 0pt{\vskip -4.5mm\line{\small\phead\ifnum
\count0=\startpage ISSN 1364-0380 (on line)
1465-3060 (printed) \hfill {\pnum\folio}\else\ifodd\count0\def\\{ }%
\ifx\theshorttitle\relax\thetitle\else\theshorttitle\fi\hfill{\pnum\folio}
\else\def\\{ and }{\pnum\folio}\hfill\ifx\theshortauthors\relax\theauthors
\else\theshortauthors\fi\fi\fi}\vss}}
\footline{\vbox to 0pt{\vglue 0mm\line{\small\pfoot\ifnum\count0=\startpage
\copyright\ \gtp\hfill\else
\gt, Volume \thevolumenumber\ (\thevolumeyear)\hfill\fi}\vss
}}
\else
\makeatletter
\def\@oddhead{{\small\lhead\ifnum\count0=\startpage ISSN 1364-0380 (on line)
1465-3060 (printed) \hfill {\lnum\number\count0}\else\ifodd\count0
\def\\{ }\ifx\theshorttitle\relax \thetitle \else\theshorttitle\fi\hfill
{\lnum\number\count0}\else\def\\{ and }{\lnum\number\count0}
\hfill\ifx\theshortauthors\relax 
\theauthors\else\theshortauthors\fi\fi\fi}}\def\@evenhead{\@oddhead}
\def\@oddfoot{\small\lfoot\ifnum\count0=\startpage\copyright\ \gtp\hfill\else
\gt, Volume \thevolumenumber\ (\thevolumeyear)\hfill\fi}
\def\@evenfoot{\@oddfoot}
\makeatother
\fi

\newwrite\gtoutfile
\long\gdef\makeheadfile{  
{\def\\{, }\def\s{ }
\immediate\openout\gtoutfile head.xxx
\immediate\write\gtoutfile{Proxy-for: \ifx\theasciiauthors\relax
\theauthors\else\theasciiauthors\fi\s<\ifx\theasciiemail\relax\theemail\else\theasciiemail\fi>}
\immediate\write\gtoutfile{\noexpand\\}
\immediate\write\gtoutfile{Authors: \ifx\theasciiauthors\relax
\theauthors\else\theasciiauthors\fi}
{\def\\{ }\immediate\write\gtoutfile{Title: \ifx\theasciititle\relax
\thetitle\else\theasciititle\fi}}
\immediate\write\gtoutfile{Subj-class: GT or SG or MG etc}
\immediate\write\gtoutfile{MSC-class: \theprimaryclass\ifx\thesecondaryclass\relax\else, \thesecondaryclass\fi}
\immediate\write\gtoutfile{Journal-ref: Geom. Topol. \thevolumenumber
(\thevolumeyear) \startpage-\finishpage}
\immediate\write\gtoutfile{Comments: Published by Geometry and Topology at}
\immediate\write\gtoutfile{\s\s http://www.maths.warwick.ac.uk/gt/GTVol\thevolumenumber/paper\thepapernumber.abs.html}
\immediate\write\gtoutfile{\noexpand\\}
\immediate\write\gtoutfile{}
\ifx\theasciiabstract\relax
\immediate\write\gtoutfile{\theabstract}\else
\immediate\write\gtoutfile{\theasciiabstract}\fi
\immediate\write\gtoutfile{}
\immediate\write\gtoutfile{\noexpand\\}
\immediate\write\gtoutfile{}
\immediate\closeout\gtoutfile}}  

\def\maketitlepage{\maketitlep\makeheadfile}
\let\maketitle\maketitlepage

\lognumber{624}
\received{11 June 2005}
\volumenumber{9}\papernumber{44}\volumeyear{2005}
\pagenumbers{1953}{1976}   
\published{17 October 2005}
\accepted{10 October 2005}
\proposed{Martin Bridson}
\seconded{Wolfgang L\"uck, Haynes Miller}

\usepackage{amsmath,amssymb}

\newtheorem{theorem}{Theorem}
\newtheorem{proposition}[theorem]{Proposition}
\newtheorem{lemma}[theorem]{Lemma} 
\newtheorem{corollary}[theorem]{Corollary}

\theoremstyle{definition} 

\newtheorem*{remark}{Remark}

\newcommand{\semi}{{:}}
\newcommand{\tr}{{\mathrm{tr}}}
\newcommand{\cc}{{\mathbb C}}
\newcommand{\ff}{{\mathbb F}}
\newcommand{\qq}{{\mathbb Q}}
\newcommand{\rr}{{\mathbb R}}
\newcommand{\zz}{{\mathbb Z}}
\newcommand{\supp}{{\mathrm{supp}}}

\newcommand\mapright[1]{\smash{\mathop{\longrightarrow}\limits^{#1}}}

\begin{document} 
\title{On finite subgroups of groups of type VF}
\author{Ian J Leary}
\address{Department of Mathematics, The Ohio State 
University\\231 West 18th Avenue, Columbus, OH 43210-1174, USA} 
\secondaddress{School of Mathematics, University of 
Southampton\\Southampton, SO17 1BJ, UK} 
\asciiaddress{Department of Mathematics, The Ohio State 
University\\231 West 18th Avenue, Columbus, OH 43210-1174, 
USA\\and\\School of Mathematics, University of 
Southampton\\Southampton, SO17 1BJ, UK} 
\email{leary@math.ohio-state.edu} 

\begin{abstract} 
For any finite group $Q$ not of prime power order, we construct 
a group $G$ that is virtually of type $F$, contains infinitely many
conjugacy classes of subgroups isomorphic to $Q$, and contains only 
finitely many conjugacy classes of other finite subgroups.  
\end{abstract}

\asciiabstract{%
For any finite group Q not of prime power order, we construct 
a group G that is virtually of type F, contains infinitely many
conjugacy classes of subgroups isomorphic to Q, and contains only 
finitely many conjugacy classes of other finite subgroups.}

\asciikeywords{Conjugacy classes, finite subgroups, groups of type F} 
\maketitle

\section{Introduction} 
A group $H$ is said to be of type $F$ if there is a finite classifying
space for~$H$, ie, if there exists a finite simplicial complex whose 
fundamental group is isomorphic to $H$ and whose universal cover is
contractible.  A group of type $F$ is necessarily torsion-free.  
It is easily seen that any finite-index subgroup of a group of type $F$
is also of type $F$.  

A group $G$ is said to be of type $VF$ if $G$ contains a finite-index
subgroup $H$ which is of type $F$, ie, if $G$ is virtually of type
$F$.  If $H$ has index $n$ in $G$, then the kernel of the action of 
$G$ on the cosets of $H$ has index at most~$n!$.  Hence any group of
type $VF$ contains a finite-index normal subgroup of type~$F$, and so
for any group $G$ of type $VF$ there is a bound on the orders of
finite subgroups of $G$.  

K\,S Brown's book `Cohomology of Groups' contains a result that
implies that a group of type $VF$ can contain only finitely many
conjugacy classes of subgroups of prime power order 
\cite[IX.13.2]{brown}.  The question of whether a group of 
type $VF$ could
ever contain infinitely many conjugacy classes of finite subgroups was
posed in \cite{wall,lueck}, and remained unanswered until Brita
Nucinkis and the author constructed examples in~\cite{vfg}.  These 
examples may be summarized as follows: 

\begin{theorem} \label{vfg}
Let $Q$ be a finite group admitting a simplicial
action on a finite contractible simplicial complex $L$ such that 
the fixed point set $L^Q$ is empty.  Then there is a group $H_L$ 
of type $F$ (depending only on $L$) and an action of $Q$ on $H_L$ 
such that the semi-direct product $H_L\semi Q$ contains infinitely 
many conjugacy classes of subgroup isomorphic to $Q$.  
\end{theorem} 

R Oliver has shown that a finite group $Q$ admits an action on a
finite contractible $L$ without a global fixed point if and only if 
$Q$ is not expressible as $p$--group-by-cyclic-by-$q$--group for any 
primes $p$~and~$q$~\cite{bobo}.  (Oliver's main result is 
the construction of actions: the proof that actions do not exist in 
the other cases is far simpler and we include it in
Section~\ref{final}.)   

The purpose of this paper is to close the gap between Brown's result 
and the construction of Theorem~\ref{vfg}.  
For any finite group $Q$ that is not of prime power order, we
construct a group $H$ of type $F$ with an action of $Q$ so that the
semi-direct product $H\semi Q$ contains infinitely many conjugacy
classes of subgroup isomorphic to $Q$, and finitely many conjugacy
classes of other finite subgroups.  As a corollary we obtain the 
following apparently stronger result.  

\begin{theorem} \label{freeprod} Let ${\cal Q}=\{Q_1,\ldots, Q_n\}$ be
a finite list of isomorphism types of finite group, such that no $Q_i$
is a group of prime power order.  There exists a group $G=G({\cal Q})$ 
of type $VF$ such that $G$ contains infinitely many conjugacy classes of
subgroup isomorphic to a finite group $Q$ if and only if $Q\in {\cal Q}$.  
\end{theorem}

In particular, it follows that a group of type $VF$ may contain
infinitely many conjugacy classes of \emph{elements} of finite order, 
although any such group can only contain finitely many conjugacy 
classes of elements of prime power order.  

Our techniques also apply to other weaker finiteness conditions.
Recall that a group $G$ is of type $FP$ over a ring $R$ if 
the trivial module $R$ for the group ring $RG$ admits a finite
resolution by finitely generated projective $RG$--modules, ie, if and
only if there is an integer $n$ and an exact sequence of $RG$--modules 
\[0\rightarrow P_n\rightarrow \cdots \rightarrow P_1\rightarrow P_0
\rightarrow R\rightarrow 0\] 
in which each $P_i$ is a finitely generated projective.  If there
exists such a sequence in which each $P_i$ is a finitely generated
free module, then $G$ is said to be of type $FL$ over $R$.  

In~\cite{vfg} Brita Nucinkis and the author proved the following.  

\begin{theorem} \label{fpq}
Let $Q$ be a finite group admitting a simplicial
action on a finite $\qq$--acyclic simplicial complex $L$ such that 
the fixed point set $L^Q$ is empty.  Then there is a virtually
torsion-free group $G=H_L\semi Q$ of type $FP$ over $\qq$ containing 
infinitely many conjugacy classes of subgroup isomorphic to $Q$.  
\end{theorem} 

R Oliver has shown that a finite group $Q$ admits such an action if
and only if $Q$ is not of the form cyclic-by-$p$--group for some prime
$p$~\cite{bobo}.  In particular, the above construction did not give 
rise to any groups of type $FP$ over $\qq$ containing infinitely many
conjugacy classes of \emph{elements} of finite order.  The question of
whether such groups can exist was posed by H Bass in~\cite{bass,wall}.  
One reason why this question is of interest is that if $G$ contains 
infinitely many conjugacy classes of elements of finite order, then 
the Grothendieck group $K_0(\qq G)$ may be shown to have infinite
rank.  (We give a proof of this fact below in Theorem~\ref{bass}.)  

Any group of type $F$ is of type $FP$ over any ring $R$, and a 
group $G$ of type $VF$ is of type $FP$ over any ring $R$ in which the 
orders of all finite subgroups of $G$ are units.  In particular, 
every group of type $VF$ is of type $FP$ over $\qq$.  
It follows that examples coming from Theorem~\ref{freeprod} 
may be used to answer Bass's question.  
By Brown's result, groups of type $VF$ necessarily 
contain only finitely many conjugacy classes of elements of prime 
power order.  This is not the case for groups of type $FP$ over $\qq$, 
and in fact for any non-trivial finite group $Q$ we construct a 
group of type $FP$ over $\qq$ containing infinitely many conjugacy 
classes of subgroup isomorphic to $Q$, and finitely many conjugacy 
classes of other finite subgroup.  The following is a corollary
of our result.  

\begin{theorem} \label{freepq} Let ${\cal Q}=\{Q_1,\ldots, Q_n\}$ be
a finite list of isomorphism types of non-trivial finite groups.  
There exists a virtually torsion-free group $G=G({\cal Q})$ of type 
$FP$ over $\qq$ 
such that $G$ contains infinitely many conjugacy classes of
subgroup isomorphic to a finite group $Q$ if and only if $Q\in {\cal Q}$.  
\end{theorem}

The groups $H_L$ appearing in the statements of Theorems
\ref{vfg}~and~\ref{fpq} are the groups introduced by M Bestvina 
and N Brady, who used them to solve a number of open problems 
concerning homological finiteness conditions~\cite{bb}.  In 
particular, in the case when $L$ is a finite acyclic
complex that is not contractible, they showed that the group $H_L$ is 
of type $FL$ over $\zz$ but is not finitely presented.  The main idea
in~\cite{vfg} was to allow a finite group $Q$ to act on the complex
$L$, and hence on the group $H_L$.  

The main idea in this paper is to consider Bestvina--Brady groups $H_L$
for \emph{infinite} complexes $L$.  
If $Q$ is any finite group not of prime power
order, then there exists a complex $L$ with a $\zz\times Q$--action such that 
\begin{enumerate} 
\item $L$ is contractible; 

\item $\zz\times Q$ acts cocompactly on $L$; 

\item all cell stabilizers are finite; 

\item $\{0\}\times Q$ fixes no point of $L$.  

\end{enumerate} 
The first three properties together imply that the semi-direct product 
$H_L\semi\zz$ is of type $F$, and the fourth property implies that the 
semi-direct product $H_L\semi(\zz\times Q)$ contains infinitely many
conjugacy classes of subgroup isomorphic to $Q$.  
A construction for $L$ as above in the case when $Q$ is cyclic was 
given by Conner and Floyd~\cite{cofl}.  In Section~\ref{final} we give a 
construction for arbitrary finite $Q$ which was shown to us by Bob Oliver.  

A similar (but simpler) construction involving an infinite
$\qq$--acyclic complex $L$ is used in proving our theorem concerning 
groups of type $FP$ over $\qq$.  

In the final section of the paper we discuss some further finiteness 
properties of the groups that we construct.  
We show that the groups are residually finite, although 
we are unable to decide whether they are linear.  We also show that 
each of the groups used in the proofs of Theorems
\ref{freeprod}~and~\ref{freepq} occurs as the kernel of a map to 
$\zz$ from a group that acts cocompactly with finite stabilizers on a 
CAT(0) cube complex.  

The work in this paper builds on the author's joint work with Brita
Nucinkis and uses theorems concerning actions of finite groups which 
the author learned from Bob Oliver.  The author gratefully
acknowledges their contributions to this work.  Some of this work was 
done at Paris~13 and at the ETH, Z\"urich.  The author thanks these 
institutions for their hospitality.  

The author was partially supported by NSF grant DMS-0505471.

\section{Bestvina--Brady groups} \label{bebr}
In this section we define the Bestvina--Brady group $H_L$ 
associated to a flag complex $L$, and we check that some of the 
results in~\cite{bb,vfg} extend to the case when $L$ is an infinite
flag complex.  

A flag complex, $L$, is a simplicial complex which contains as many 
higher dimensional simplices as possible, given its 1--skeleton.  
In other words, whenever the complete graph on a finite subset of 
the vertex set of $L$ is contained in the 1--skeleton of $L$, then 
there is a simplex of $L$ with that set of vertices.  The 
realisation of any poset is a flag complex (since a subset is totally
ordered if any two of its members are comparable).  In particular, 
the barycentric subdivision of any simplicial complex is a flag
complex.  

Given a flag complex $L$, the associated right-angled Artin group 
$G_L$ is the group with generators the vertices of $L$ subject only 
to the relations that the ends of each edge commute.  There is a 
model for the classifying space $BG_L$ with one $n$--dimensional
cubical cell for each $(n-1)$--simplex of $L$ (including one vertex 
corresponding to the empty simplex in $L$).  Let $X_L$ denote the 
universal cover of this space.  Cells of $X_L$ are 
$n$--cubes of the form $(g,v_1,\ldots,v_n)$
where $(v_1,\ldots,v_n)$ is an $n-1$--simplex of $L$ and $g$ is an 
element of $G$.  The $i$th pair of opposite faces of this $n$--cube
consists of the cubes $(g,v_1,\ldots,\hat v_i,\ldots, v_n)$ and 
$(gv_i,v_1,\ldots,\hat v_i,\ldots,v_n)$, where $gv_i$ is the product 
of two elements of $G_L$, and as usual $\hat v_i$ means `omit $v_i$'.  
The action of $G_L$ is given by $g'(g,v_1,\ldots,v_n)=
(g'g,v_1,\ldots,v_n)$.  If $\sigma= (v_1,\ldots,v_n)$ is a simplex of 
$L$, we will write $(g,\sigma)$ in place of $(g,v_1\ldots,v_n)$.  In 
particular, we will write $(g)$ for a vertex of $X_L$.  

The space $X_L$ admits the structure of a 
CAT(0) cubical complex:  there is a geodesic CAT(0) metric on $X_L$ in 
which each cubical cell is isometric to a standard Euclidean unit
cube, and the action of $G_L$ is by isometries of this metric.  In 
the case when $L$ is infinite, $X_L$ is not locally finite, and the 
metric topology on $X_L$ is not the same as the CW--topology, 
but this will not cause any difficulties.  

Suppose now that $f\co L\rightarrow L'$ is a simplicial map.  Then $f$
defines a group homomorphism $f_*\co G_L\rightarrow G_{L'}$, which takes the
generator $v$ to the generator $f(v)$, and $f$ defines a
piecewise-linear continuous map $f_!\co X_L\rightarrow X_{L'}$, which 
takes the vertex $(g)$ to the vertex $(f(g))$, and extends linearly
across each cube.  The map $f_!$ is $G_L$--equivariant, where $f_*$ is
used to define the $G_L$--action on $X_{L'}$, and so induces a map 
from $X_L/G_L$ to $X_{L'}/G_{L'}$, which is an explicit construction
for the map $B(f_*)\co BG_L\rightarrow BG_{L'}$.  If $f$ 
is an injective simplicial map, then $f_*$ is an injective group
homomorphism and $f_!$ is an isometric embedding.  

Two special cases of this construction are of interest to us. 
Firstly, any group $\Gamma$ of automorphisms of $L$ gives rise to 
a group of automorphisms of $G_L$ and to a group of cellular 
automorphisms of $X_L/G_L$.  Since the unique vertex in $X_L/G_L$ 
is fixed by $\Gamma$, the group of all lifts of elements of $\Gamma$
to the covering space $X_L\rightarrow X_L/G_L$ is the semi-direct 
product $G_L\semi\Gamma$, where $\Gamma$ acts on $G_L$ via the 
action described above.  

Secondly, let $*$ denote a 1--point simplicial complex.  For this 
choice of simplicial complex, $G_*$ is infinite cyclic, and $X_*$ is
the real line triangulated with one orbit of vertices and one orbit of
edges.  For any $L$, there is a unique map $f_L\co L\rightarrow *$, and the 
Bestvina--Brady group $H_L$ is defined to be the kernel of
$f_{L*}\co G_L\rightarrow G_*$.  The map $f_{L!}\co X_L\rightarrow X_*\cong \rr$ 
may be considered as defining a `height function' on $X_*$.
Identifying the integers $\zz\subseteq \rr$ with the vertex set in $X_*$, 
one sees that $f_{L!}$ sends each vertex of $X_L$ to an integer, and 
that each cube of $X_L$ has a unique minimal and maximal vertex for 
this height function.  For the cube $C$, we shall write $\min(C)$ and 
$\max(C)$ respectively for its minimal and maximal vertices.  
Any simplicial map $f\co L\rightarrow L'$ fits in
to a commutative triangle with $f_L\co L\rightarrow *$ and
$f_{L'}\co L'\rightarrow *$, and hence one obtains an induced map
$f_*\co H_L\rightarrow H_{L'}$.  In particular, if $\Gamma$ is a
group of simplicial automorphisms of $L$, then the semi-direct product 
$H_L\semi\Gamma$ is defined and is equal to the kernel of the 
composite $G_L\semi\Gamma\rightarrow G_*\times \Gamma\rightarrow
G_*$.  

The work of Bestvina and Brady~\cite{bb} relies on a study of 
the height function $f\co X_L\rightarrow X_*=\rr$.  We recall part 
of this, and check that it applies to the case when $L$ is 
infinite (which was not considered in~\cite{bb}).  

Pick a point $c$ in the interior of an edge of $X_*$, and 
define $Y=Y_L= f^{-1}(c)\subseteq X_L$.  
(The point $c$ will remain fixed for 
the remainder of this section, but will be suppressed from the 
notation.)  Give $Y$ the structure of a polyhedral 
CW--complex by taking as cells the sets $C\cap Y$ where $C$ is a 
cube of $X_L$.  Note that the CW--structure on $Y$ gives rise to the
same topology as the subspace topology coming from the CW--topology on
$X$.  

Now let $C$ be a cube in $X_L$ whose highest
vertex is $v_1$.  Define a subset $C_c$ of $C$ by 
$$C_c = C \cap f^{-1}([c,\infty)) = C \cap f^{-1}([c,f(v_1)]).$$ 
Similarly, if the lowest vertex of $C$ is $v_0$, define a subset 
$C^c$ by 
$$C^c = C \cap f^{-1}((-\infty,c]) = C \cap f^{-1}([f(v_0),c]).$$ 
If $C= (g,\sigma)$ for some simplex $\sigma\in L$, then the link 
of $v_1$ in $C$ is homeomorphic to $\sigma$.  It follows that 
if $f(v_1)>c$, then $C_c$ is homeomorphic to the cone on $L$ with 
vertex $v_1$.  If $f(v_1)<c$, then $C_c$ is empty.  Similarly, if 
$f(v_0)<c$ then $C^c$ is empty, and otherwise $C^c$ is homeomorphic to
the cone on $\sigma$.  Now for $v$ a vertex of $X_L$, define $F(v)$ to 
be either 
\[F(v) = \begin{cases}\bigcup_{v=\max(C)}\,\, C_c & f(v)>c\cr 
 \bigcup_{v=\min(C)}\,\, C^c & f(v)<c\end{cases}
\] 
For each $v$, one may show that $F(v)$ is homeomorphic to the cone on 
$L$ with vertex $v$.  (Here, as usual, we are using the CW--topology on
both $F(v)$ and the cone on $L$.)  
Now for $a,b\in X_*=\rr$ with $a<c<b$, define a
subspace $Y_{[a,b]}$ of $X_L$ by 
\[Y_{[a,b]} = Y\cup \bigcup_{a\leq f(v)\leq b}F(v).\]
Each $Y_{[a,b]}$ is a CW--complex, with cells the truncated cubes 
$C_c$, $C^c$ and $C\cap Y$ for each cube $C$ of $X_L$, and if 
$\alpha\leq a<c<b\leq \beta$, then $Y_{[a,b]}$ is a subcomplex of 
$Y_{[\alpha,\beta]}$.  As $a$ decreases (resp.\ $b$ increases) 
the complex $Y_{[a,b]}$ only changes as $a$ (resp.\ $b$) passes
through an integer.  For each $a<c<b$, one has that $Y_{[a-1,b+1]}$ 
is homeomorphic to $Y_{[a,b]}$ with a family of subspaces
homeomorphic to $L$ coned off.  (There is one such cone for 
each vertex in $Y_{[a-1,b+1]}-Y_{[a,b]}$.)  
Thus one obtains the following lemma 
and corollary due to Bestvina--Brady~\cite{bb}, for any 
simplicial complex $L$.  

\begin{lemma} If $L$ is contractible, then for any $a<c<b$, the 
inclusion of $Y$ in $Y_{[a,b]}$ is a homotopy equivalence.  If 
$L$ is $R$--acyclic for some ring $R$, then for any $a<c<b$, the 
inclusion of $Y$ in $Y_{[a,b]}$ induces an isomorphism of
$R$--homology.  
\end{lemma}

\begin{corollary} \label{besbra} 
If $L$ is contractible, then $Y$ is contractible.
If $L$ is $R$--acyclic, then $Y$ is $R$--acyclic.  
\end{corollary} 

\begin{proof} We know that $X_L$ is contractible, and the lemma
implies that the inclusion $Y\rightarrow X_L$ is a homotopy
equivalence if $L$ is contractible and is an $R$--homology isomorphism
if $L$ is $R$--acyclic.  
\end{proof}  

\begin{theorem}\label{type} 
Suppose that $\Gamma$ acts freely cocompactly on a
simplicial complex~$L$.  If $L$ is contractible, then
$H_L\semi\Gamma$ is type $F$.  If $L$ is $R$--acyclic, then
$H_L\semi\Gamma$ is type $FL$ over $R$.  
\end{theorem}  

\begin{proof} It follows 
from Corollary~\ref{besbra} that $Y$ is contractible or
$R$--acyclic whenever $L$ is.  Thus it suffices to show that
$H_L\semi\Gamma$ acts freely cocompactly on $Y$.  To see this, first 
note that $G_L\semi\Gamma$ has only finitely many orbits of cells in 
its action on $X_L$.  If $C$ is an $n$--cube of $X_L$ with top vertex 
$v$, then $C\cap Y$ is non-empty if and only if $c<f(v)< c+n$.  It 
follows that each $G_L\semi\Gamma$--orbit of $n$--cubes in $X_L$ gives 
rise to exactly $n$ $H_L\semi\Gamma$--orbits of $(n-1)$--cells in $Y$.  
\end{proof} 

It remains to study the conjugacy classes of finite subgroups of 
groups of the form $H_L\semi\Gamma$ and $G_L\semi\Gamma$.  In fact 
it is no more difficult to study conjugacy classes of subgroups 
$Q'$ such that $Q'\cap G_L= \{1\}$.  Consider first 
the collection of subgroups $\Gamma'$ of $G_L\semi\Gamma$ 
which map isomorphically to $G_L\semi\Gamma/G_L\cong \Gamma$.  
The action of $\Gamma$ on $X_L/G_L$ fixes the unique vertex.  It 
follows that each such $\Gamma'$ fixes some vertex $v$ of $X_L$.  
Since the vertices form a single orbit, it follows that all such 
$\Gamma'$ are conjugate in $G_L\semi\Gamma$.  

\begin{proposition} Let $\Gamma$ act on $L$, let $Q\leq \Gamma$, 
and let $Q'$ be any subgroup of $G_L\semi\Gamma$ that maps 
isomorphically to $Q\leq \Gamma= G_L\semi\Gamma/G_L$.    
If $L^Q=\emptyset$, then $Q'$ fixes a unique vertex in $X_L$.  
If $L^Q$ contains the barycentre of an $m$--simplex, and $Q'$ fixes a
vertex $(g)\in X_L$ of height $f(g)=a$, then $Q'$ also fixes a 
vertex of height $a+(m+1)n$ for each integer $n$.  
\end{proposition} 

\begin{remark} Since we are not assuming that the action of $\Gamma$ 
on $L$ makes $L$ into a $\Gamma$--CW--complex, it is not necessarily 
the case that $L^Q$ is a subcomplex of $L$.  However there can be a
point of $L^Q$ in the interior of the simplex $\sigma$ only if 
$q\sigma=\sigma$ for all $q\in Q$.  In this case the barycentre of 
$\sigma$ is a point fixed by $Q$.  
\end{remark} 

\begin{proof} For the first time, we shall make use of the CAT(0)
metric on $X_L$.  Suppose that $Q'$ fixes two distinct vertices 
$(g),(h)$ of $X_L$.  Since geodesics in a CAT(0) metric space are 
unique, it follows that the geodesic arc from $(g)$ to $(h)$ is 
also fixed by $Q'$.  The start of this arc is a straight line 
passing from $(g)$ into the interior of $C$, an $n$--cube of $X_L$ 
for some $n>0$, which must be preserved (setwise) by $Q'$.  If 
$C=(g',v_1,\ldots,v_n)$, then it follows that the $(n-1)$--simplex 
$(v_1,\ldots,v_n)$ in $L$ is (setwise) preserved by $Q$, and hence 
that $L^Q\neq \emptyset$.  

For the second statement, suppose that $(g)$ is fixed by $Q'$, and
that the $m$--simplex $(v_0,\ldots,v_m)$ in $L$ is setwise fixed by
$Q$.  Then the long diagonal from $(g)$ to $(gv_0v_1\cdots v_m)$ in 
the $(m+1)$--cube $(g,v_0,\ldots,v_m)$ is an arc fixed by $Q'$, which
connects two vertices whose heights differ by $m+1$.  It follows 
that for any $n$, the vertex $g(v_0v_1\cdots v_m)^n$ is fixed by
$Q'$.  
\end{proof} 

\begin{theorem} \label{conjclass}
Let $\Gamma$ act on $L$, and let $Q\leq \Gamma$. 
If $L^Q=\emptyset$, then there are infinitely many conjugacy classes
of subgroups $Q'$ of $H_L\semi\Gamma$ whose members map isomorphically 
to conjugates of $Q$ in $\Gamma$.  
If $L^Q$ contains the barycentre of an $m$--simplex, 
then there are at most $m+1$ conjugacy classes of such $Q'$ in 
$H_L\semi\Gamma$.  In particular, if $L^Q$ contains a vertex of $L$, 
then any two such subgroups are conjugate.  
\end{theorem} 

\begin{proof} 
We know that any such $Q'$ fixes a vertex of $X_L$ and 
that every vertex is fixed by some such $Q'$.  
In the case when $L^Q=\emptyset$, each $Q'$ fixes exactly one vertex 
of $X_L$.  Since vertices of different heights are in different orbits 
for the action of $H_L\semi\Gamma$, it follows that in this case there 
are infinitely many conjugacy classes of such $Q'$.  

In general, $H_L$ acts transitively on the vertices of fixed height.  
If $L^Q$ contains the barycentre of an $m$--simplex, and $Q'$ fixes a 
vertex of height $a$, then $Q'$ also fixes a vertex of height
$a+(m+1)n$ for each $n$.  Hence given any $(m+2)$ subgroups of
$H_L\semi\Gamma$ which map isomorphically to $Q$ or one of its
conjugates, some pair $Q'$, $Q''$ of these subgroups must fix 
vertices of the same height.  Let $\Gamma'\geq Q'$ and $\Gamma''\geq
Q''$ be the stabilizers of these vertices, which map isomorphically to 
$\Gamma$.  The subgroups $\Gamma'$ and $\Gamma''$ are conjugate by 
an element of $H_L$.  Hence it follows that $Q'$ and $Q''$ are 
conjugate by some element of $H_L\semi\Gamma$.  
\end{proof}

\section{Group actions} 

\label{final}
Here we construct the actions of finite groups $Q$ and direct products
of the form $\zz\times Q$ on finite-dimensional simplicial complex
that are needed in order to apply the constructions of the previous
section.  The first two propositions are included to show why actions 
of finite groups alone cannot give all the examples that we need.  

\begin{proposition} Suppose that $Q$ is a finite group with 
normal subgroups $P\leq P'$, so that $P$ and $Q/P'$ are 
groups of prime power order and $P'/P$ is cyclic.  For any 
action of $Q$ on a finite contractible complex $L$, the fixed point
set $L^Q$ is non-empty.  
\end{proposition}

\begin{proof} Let $p$ and $q$ be the primes (not necessarily distinct) 
so that $|P|$ is a power of $p$ and $|Q:P|$ is a power of $q$.  Let 
$C$ denote $P'/P$, and let $Q'$ denote $Q/P'$.  

By P\,A Smith theory~\cite[VII.10.5]{brown}, 
the fixed point set $L'=L^P$ has the same mod-$p$
homology as a point, and hence has Euler characteristic equal to~1.  
By character theory, it follows that the Euler characteristic of 
$L''=L^{P'}= {L'}^C$ is equal to~1.  By counting lengths of orbits of
cells, one sees that $L^Q= {L''}^{Q'}$ has Euler characteristic congruent
to~1 modulo~$q$.  This implies that $L^Q$ is not empty.  
\end{proof} 

The above proof also gives:

\begin{proposition} Let $Q$ be a finite group with a normal cyclic
subgroup $P'$ so that $Q/P'$ is a group of prime power order.  For 
any action of $Q$ on a finite complex $L'$ with Euler characteristic 
$\chi(L')=1$, the fixed point set ${L'}^Q$ is non-empty. 
\end{proposition}

The actions on $\qq$--acyclic spaces that we shall need will all come
from Theorem~\ref{qacyc}.  In the proof of we shall need
Lemma~\ref{wall} concerning Wall's finiteness obstruction.  

Suppose that $G$ is a group of type $FP$ over the ring $R$, and that 
$$0\rightarrow P_n\rightarrow \cdots \rightarrow P_1\rightarrow P_0
\rightarrow R\rightarrow 0$$ 
is a resolution of $R$ over $RG$ by finitely generated projectives.
As usual, let $K_0(RG)$ denote the Grothendieck group of
finitely-generated projective $RG$--modules.  The Wall obstruction 
or Euler characteristic of $G$ over $R$ is the element of $K_0(RG)$ 
given by the alternating sum 
$$w(R,G)= \sum_i (-1)^i[P_i]$$
and is independent of the choice of resolution \cite[I.7]{rosenberg}.  

\begin{lemma} Let $Q$ be a finite group.  The group $G=\zz\times Q$
is $FP$ over $\qq$, and the Wall obstruction for this group is zero.  
\label{wall}
\end{lemma} 

\begin{proof} 
Let the group $G=\zz\times Q$ act on the real line via the 
projection $G\rightarrow \zz$.  There is a $G$--equivariant
triangulation of the line with one orbit of 0--cells of type 
$G/Q$ and one orbit of 1--cells, also of type $G/Q$.  The 
cellular chain complex for this space gives a projective 
resolution for $\qq$ over $\qq G$ of length one: 
$$0\rightarrow 
\qq G/Q \rightarrow \qq G/Q \rightarrow \qq \rightarrow 0,$$
in which the modules in degrees 0 and 1 are isomorphic to each other.  
\end{proof} 

\begin{theorem} Let $Q$ be a finite group, and let $\cal F$ be a 
non-empty family of subgroups of $Q$ which is closed under conjugation and
inclusion.   There is a 3--dimensional $\qq$--acyclic simplicial complex 
$L$ admitting a cocompact action of $\Gamma=\zz\times Q$ so that all cell
stabilizers are finite and so that $P\leq Q$ fixes some point of $L$
if and only if $P\in \cal F$.  \label{qacyc}
\end{theorem} 

\begin{proof} Let $\Delta$ be a finite set with a $Q$--action, such
that $\Delta^P\neq \emptyset$ if and only if $P\in \cal F$, and let 
$Z=Q*Q*\Delta$ be the join of two copies of $Q$ and one copy of
$\Delta$, with the diagonal action of $Q$.  This $Z$ is a
2--dimensional simply-connected $Q$--space, with the property that the
$Q$--action is free except on the 0--skeleton.  Let $\zz$ act on $\rr$
in the usual way, and let $L_0$ be the product $\rr\times Z$ with the 
product action of $\Gamma=\zz\times Q$.  Now let $L_1$ be the 2--skeleton of 
$L_0$.  The cells of $L_1$ in non-free orbits form a copy of
$\rr\times \Delta$ with the product action of $\Gamma$.  Let 
$C_*$ be the rational chain complex for $L_1$.  Since $L_1$ is
1--connected, $C_*$ forms the start of a projective resolution for 
$\qq$ over $\qq\Gamma$.  As $\qq\Gamma$--modules, 
$C_2$ is free and each of $C_1$ and $C_0$ is the direct
sum of a free module and a copy of $\qq[\zz\times\Delta]$.  
Hence the element of $K_0(\qq\Gamma)$ defined by the alternating sum 
$[C_2]-[C_1]+[C_0]$ is in the subgroup of $K_0(\qq\Gamma)$ generated 
by the free module.  Since we know by Lemma~\ref{wall} that the 
Wall obstruction for $\Gamma$ over $\qq$ is zero, it follows that 
$H_2(C_*)$ is a stably-free $\qq\Gamma$--module.  Make $L_2$ by
attaching finitely many free $\Gamma$--orbits of 2--spheres to $L_1$
in such a way that $H_2(L_2;\qq)$ is a free $\qq\Gamma$--module.  
Let $c_1,\ldots,c_k$ be cycles in $C_2(L_2,\qq)$ representing a
$\qq\Gamma$--basis for $H_2(L_2;\qq)$, and pick a large integer $M$ 
so that each $M.c_i$ is an integral cycle.  Since $L_2$ is
1--connected, each $M.c_i$ is realized by the image of the fundamental 
class for $S^2$ under some map $f_i\co S^2\rightarrow L_2$.  
Now define $L_3$ by attaching free $\Gamma$--orbits of
3--balls to $L_2$, using the $f_i$ as attaching maps for orbit
representatives.  This $L_3$ has all of the required properties,
except that it has been constructed as a $\Gamma$--CW--complex rather
than as a $\Gamma$--simplicial complex.  By the simplicial
approximation theorem, we can construct a 3--dimensional
$\Gamma$--simplicial complex $L$ together with an equivariant 
homotopy equivalence $L\rightarrow L_3$.  
\end{proof} 

Before stating and proving Theorem~\ref{contr}, which will provide all
the actions on contractible spaces that we shall need, we begin by
establishing some notation, and proving a lemma concerning equivariant
self-maps of spheres in linear representations.  Lemma~\ref{sphe} and 
Theorem~\ref{contr} were shown to the author by Bob Oliver.  

Let $S$ denote the unit sphere in
$\cc^n$, so that $S$ is a sphere of odd dimension.  For $x\in S$, let 
$T_xS$ be the tangent space to $S$ at $x$, and let $B_x$ be the 
closed unit ball in $T_xS$, with boundary $\partial B_x$.  For
$\epsilon>0$, let $e_{\epsilon,x}\co B_x\rightarrow S$ denote the scalar
multiple of the exponential map such that the image of $B_x$ is a ball 
of radius $\epsilon$ in $S$.  In the case when $\epsilon=\pi$, this
map sends the whole of $\partial B_x$ to the point $-x$.  The cases 
of interest to us include the case $\epsilon =\pi$ and the case when 
$\epsilon$ is small.  Suppose that a finite group $P$ acts linearly 
on $S$, fixing the point $x$.  This induces a $P$--action on $T_xS$, 
and the exponential map $e_{\epsilon,x}$ is $P$--equivariant in the 
sense that $e_{\epsilon,x}(gv)= ge_{\epsilon,x}(v)$ for all $v\in B_x$ 
and all $g\in P$.  

Each of the self-maps of spheres that we shall construct will have the
property that it is equal to the identity except on a number of small 
balls.  For such a map $f\co S\rightarrow S$, its support, $\supp(f)$, is 
defined to be the closure of the set of points $x\in S$ so that
$f(x)\neq x$.  Given another such map $f'\co S\rightarrow S$ with 
$\supp(f)\cap \supp(f')=\emptyset$, the map $f\coprod f'$ is defined
by 
$$f\textstyle{\coprod} f' (x) = \begin{cases} 
f(x) & \,\, x\in \supp(f)\cr 
f'(x) & \,\, x\in \supp(f')\cr 
x&\,\, x\notin \supp(f)\cup \supp(f').
\end{cases} 
$$
Suppose that a group $Q$ acts on $S$.  For $f\co S\rightarrow S$ a
self-map of $S$ and $g\in Q$, define $g*f(s) = g(f(g^{-1}(s)))$.  
The support of $g*f$ is equal to $g.\supp(f)$.  

For $x\in S$, let $r\co (B_x,\partial B_x)\rightarrow (B_x,\partial B_x)$ be
any map of degree $-1$, for example a reflection in a hyperplane
through $0$ in $B_x$.  Define $\tilde\phi_x,\tilde\psi_x: 
B_x\rightarrow S$ by 
\begin{eqnarray*}
\tilde\phi_x(v) = \begin{cases} 
-e_{\pi,x}(2v)& |v| \leq 1/2\cr 
(|v|-1/2)v & |v| \geq 1/2\cr 
\end{cases} \cr
\tilde\psi_x(v) = \begin{cases} 
-e_{\pi,x}(r(2v))& |v| \leq 1/2\cr 
(|v|-1/2)v & |v| \geq 1/2\cr 
\end{cases} \cr
\end{eqnarray*}
and define self-maps $\phi_{\epsilon,x}$ and $\psi_{\epsilon,x}$ to be
the identity outside of the image of $e_{\epsilon,x}$ and equal to 
$\tilde\phi_x\circ e_{\epsilon,x}^{-1}$ and 
$\tilde\psi_x\circ e_{\epsilon,x}^{-1}$ respectively on their
supports.  If $f$ is a self-map of $S$ of degree $n$ whose support is 
disjoint from the $\epsilon$--ball around $x$, then $f\coprod
\phi_{\epsilon,x}$ is a self-map of degree $n+1$ and $f\coprod
\psi_{\epsilon,x}$ is a self-map of degree $n-1$.  

Suppose that a finite group $Q$ acts linearly on $S$, so that the
distance between any two points of the orbit $Q.x$ is greater than
$2\epsilon$.  If $g\in Q$, then $g*\phi_{\epsilon,x}$ and 
$\phi_{\epsilon,g.x}$
are equal.  In particular, if $g$ is an element of $Q_x$, 
the stabilizer of the point $x$, then the equation 
$g*\phi_{\epsilon,x} = \phi_{\epsilon,x}$ holds.  
Since the definition of $\psi$ involved an 
arbitrary choice of function $r$, there is no corresponding
equivariance property for the $\psi$ self-maps.  However, the map 
$g*\psi_{\epsilon,x}$ is a self-map whose support is the
$\epsilon$--ball in $S$ centred at $g.x$, and if $f$ is any self-map 
of $S$ whose support is disjoint from this ball, the coproduct 
$f\coprod g*\psi_{\epsilon,x}$ is a self-map whose degree is one less 
than that of $f$.  

For any $x\in S$, define 
$$Q.\phi_{\epsilon,x}= \coprod_{g\in Q/Q_x} g*\phi_{\epsilon,x},$$ 
for any sufficiently small $\epsilon$, where the sum runs over a 
transversal to $Q_x$ in $Q$.  For $x$ in a free $Q$--orbit, define 
$$Q.\psi_{\epsilon,x} = \coprod_{g\in Q} g*\psi_{\epsilon,x},$$ 
for small $\epsilon$.  Each of these maps is $Q$--equivariant.

\begin{lemma} Let $S$ be the unit sphere in a complex representation
  of the finite group $Q$, and suppose that $S$ contains points in 
$Q$--orbits of coprime lengths.  Then $S$ admits a $Q$--equivariant 
self-map of degree zero.  \label{sphe}
\end{lemma}

\begin{proof} 
Without loss of generality, we may suppose that $Q$ acts faithfully on
$S$.  The action of the unit circle in $\cc$ on $S$ commutes with the 
$Q$--action, and so whenever $S$ contains a $Q$--orbit of a given
length, $S$ contains infinitely many $Q$--orbits of that length.  Pick
points  $x_1,\ldots,x_m$ 
in distinct $Q$--orbits, such that the sum of the lengths of the orbits
is congruent to $-1$ modulo $|Q|$, ie, so that there exists $n$ with 
$$|Q|n = 1+ \sum_{i=1}^m |Q.x_i|.$$ 
Now pick $y_1,\ldots, y_n$ in distinct free $Q$--orbits.  Choose $\epsilon$
sufficiently small that any two points in any of these orbits are 
separated by more than $2\epsilon$.  The coproduct 
$$f= \coprod_{i=1}^mQ.\phi_{\epsilon,x_i} {\textstyle{\coprod}}
\coprod_{j=1}^n Q.\psi_{\epsilon,y_j}$$
is the required degree zero map.  
\end{proof} 

\begin{theorem}\label{contr} Let $Q$ be a finite group not of prime
power order.  Then there exists a finite-dimensional contractible 
simplicial complex $L$ with a cocompact action of $\zz\times Q$ such 
that all stabilizers are finite and such that $L^Q=\emptyset$.
Furthermore, $L$ may be chosen in such a way that $L^P\neq \emptyset$ 
for $P$ any proper subgroup of $Q$.  
\end{theorem} 

\begin{proof} Let $S$ be the unit sphere in the `reduced regular
  representation of $Q$', ie, the regular representation $\cc Q$
  minus the trivial representation.  This $S$ has the property 
that $S^Q=\emptyset$ but $S^P\neq \emptyset$ for any proper subgroup 
$P<Q$.  Since $Q$ is not of prime power order, $S$ satisfies the
  hypotheses of Lemma~\ref{sphe}, and so there exists a
  $Q$--equivariant map $f\co S\rightarrow S$ of degree zero.  

Take a $Q$--equivariant triangulation of the space $I\times S$, where 
$Q$ acts trivially on the interval $I$.  By the simplicial
approximation theorem, there is an integer $n\geq 0$ and a simplicial 
map $f'\co \{1\}\times S^{(n)}\rightarrow \{0\}\times S$ which is
equivariantly homotopic to $f$.  Now let $M$ be the $n$th barycentric
subdivision of $I\times S$ relative to $\{0\}\times S$.  This is a 
copy of $I\times S$, with the original triangulation on the subspace 
$\{0\}\times S$ and the $n$th barycentric subdivision of this
triangulation on $\{1\}\times S$.  Construct $L$ from the direct
product $\zz\times M$ by identifying $(m,1,s)$ with $(m+1,0,f'(s))$
for each $s\in S$ and $m\in \zz$.  This space $L$ is a triangulation
of the doubly
infinite mapping telescope of the map $f'\co S\rightarrow S$.  The fact
that $f'$ has degree zero implies that $L$ is contractible.  
\end{proof} 

One difference between Theorem~\ref{qacyc} and Theorem~\ref{contr} is
that the the dimension of the space constructed in Theorem~\ref{contr} 
varies with $Q$.  The final results in this section show that this 
difference cannot be avoided.  

\begin{lemma} Let $Q$ be the special linear group $SL_n(\ff_p)$ over the
  field of $p$ elements.  Let $e_1,\ldots,e_n$ be the standard basis 
for the vector space $\ff_p^n$.  Define elements
$\tau_1,\ldots,\tau_n\in Q$ by 
$$\tau_i(e_j) = \begin{cases} 
e_j & i\neq j\cr 
e_i+ e_{i+1} & i=j<n\cr 
e_n+e_1 & i=j=n.\cr 
\end{cases}
$$
The elements $\tau_1,\ldots,\tau_n$ generate $Q$, and any proper
subset of them generates a subgroup of order a power of $p$.  
\label{gens} 
\end{lemma} 

\begin{proof} 
Let $\theta$ be the cyclic permutation of the $n$ standard basis
elements, so that $\theta(e_i) = e_{i+1}$ for $i<n$ and
$\theta(e_n)=\theta_1$.  The action of $\theta$ on $Q$ by conjugation 
induces a cyclic permutation of the $\tau_i$.  

The elements $\tau_1,\ldots,\tau_{n-1}$ generate the upper triangular
matrices, which form a Sylow $p$--subgroup of $Q$.  This group contains 
each of the elementary matrices $E_{i,j}$ for $i<j$, defined by 
$$E_{i,j}(e_k) = \begin{cases} 
e_k & k\neq i\cr 
e_k + e_j & k = i.\cr
\end{cases}$$ 
Conjugation by powers of $\theta$ induces a transitive permutation of
the size $n-1$ subsets of $\tau_1,\ldots,\tau_n$.  Hence one sees that 
each such set generates a Sylow $p$--subgroup of $Q$.  

It is well-known that 
the elementary matrices $E_{i,j}$ for all $i\neq j$ form a generating
set for $Q$.  Each elementary matrix may be expressed as the conjugate
of an upper triangular elementary matrix by some power of $\theta$.  
It follows that the subgroup generated by $\tau_1,\ldots,\tau_n$
contains all elementary matrices and so is equal to $Q$.  
\end{proof} 

\begin{theorem} \label{slnp}
As in the previous lemma, let $Q=SL_n(\ff_p)$.
  Suppose that $L$ is contractible, or that $L$ is mod-$p$ acyclic, 
and that $Q$ acts on $L$ so that $L^Q=\emptyset$.  Then the dimension
of $L$ is at least $n-1$.  
\end{theorem}

\begin{proof} We may assume that $L$ is finite-dimensional, or there
  is nothing to prove. 
Let $L_i$ be the fixed point subspace for the action of
  $\tau_i$.  By P. A. Smith theory, the fixed point set for the 
action of a $p$--group on a finite-dimensional mod-$p$ acyclic space is
itself mod-$p$ acyclic.  From Lemma~\ref{gens} it follows that each 
intersection of at most $n-1$ of the $L_i$ is mod-$p$ acyclic, and 
that the intersection $L_1\cap\ldots \cap L_n$ is empty.  Let $X$ be
the union of the $L_i$.  The Mayer--Vietoris spectral sequence for the 
covering of $X$ by the $L_i$ with mod-$p$ coefficients is isomorphic 
to the spectral sequence for the covering of the boundary of an 
$(n-1)$--simplex by its faces.  It follows that the mod-$p$ homology of 
$X$ is isomorphic to the mod-$p$ homology of an $(n-2)$--sphere.  Hence 
$X$ cannot be a subspace of a mod-$p$ acyclic space of dimension strictly 
less than~$n-1$.  
\end{proof} 

\begin{remark} For a discrete group $G$, the minimal dimension of any
contractible simplicial complex admitting a $G$--action without a
global fixed point is an interesting invariant of $G$.  The above
theorem shows that this invariant can take arbitrarily large finite
values.  When $G$ is
a finite group of prime power order, the invariant takes the value
infinity.  Peter Kropholler has asked whether there are any other 
finitely generated groups $G$ for which the invariant takes the value
infinity.  
\end{remark} 

\section{Examples} \label{examps}

Here we combine the results of Sections \ref{bebr}~and~\ref{final}
to construct groups with strong homological finiteness properties that
contain infinitely many conjugacy classes of certain finite
subgroups.  

\begin{theorem} \label{ivf}
Let $Q$ be a finite group not of prime power order.
There is a group $H$ of type $F$ and a group $G=H\semi Q$ such that 
$G$ contains infinitely many 
conjugacy classes of subgroup isomorphic to $Q$ and finitely many
conjugacy classes of other finite subgroups.  
\end{theorem} 

\begin{proof} By Theorem~\ref{contr}, there is a contractible
finite-dimensional simplicial complex $L$ with a cocompact action of 
$\zz\times Q$ such that all stabilizers are finite, $L^Q=\emptyset$
and $L^P\neq \emptyset$ if $P<Q$.  Take a flag triangulation of $L$, 
and consider the Bestvina--Brady group $H_L$.  By Theorem~\ref{type},
the semi-direct product $H=H_L\semi \zz$ is of type $F$.  By
Theorem~\ref{conjclass} the group $G=H_L\semi(\zz\times Q)$ contains
infinitely many conjugacy classes of subgroups isomorphic to $Q$ and 
finitely many conjugacy classes of other finite subgroups.  
\end{proof} 

We can now prove Theorem~\ref{freeprod} as stated in the
introduction.  We first give a lemma concerning free products.  

\begin{lemma} \label{prodlem}
Let $G=G_1*\cdots *G_n$ be a free product of groups, 
and let $H_i$ be a finite-index normal subgroup of $G_i$.  There is a
bijection between conjugacy classes of non-trivial finite subgroups of
$G$ and the disjoint union of the sets of conjugacy classes of
non-trivial finite subgroups of the $G_i$.  The kernel of the map 
from $G$ to $\prod_i G_i/H_i$ is isomorphic to the free product of 
finitely many copies of the $H_i$ and a finitely-generated free
group.  
\end{lemma} 

\begin{proof} Take a classifying space $BG_i$ for each $G_i$, take a
star-shaped tree with $n$ edges whose central vertex has valency $n$, 
and make a classifying space $BG$ for $G$ by attaching the given
$BG_i$ to the $i$th boundary vertex of the tree.  Now consider the 
regular covering of this space $BG$ corresponding to the kernel of 
the homomorphism $G\rightarrow \prod_i G/H_i$.  This is a finite
covering.  The subspace of this covering lying above each $BG_i$ is a
finite disjoint union of copies of $BH_i$, and the subspace lying
above the tree is a finite disjoint union of trees.  Hence the whole 
space, which is a classifying space for the kernel, consists of a
finite number of copies of the $BH_i$'s, connected together by a 
finite number of trees.  The fundamental group of such a space is the 
free product of finitely many copies of the $H_i$ and a
finitely-generated free group.  

For the claimed result concerning conjugacy classes of finite
subgroups, we consider the tree obtained from the given expression for
$G$ as a free product.  One way to construct this tree is by
considering the universal covering space of the model for $BG$ given
above.  This consists of copies of the $EG_i$'s, connected together by
trees.  Now contract each copy of $EG_i$ to a point.  The resulting
$G$--space is contractible (since replacing $EG_i$ by a single point
does not change its homotopy type) and is 1--dimensional.  It is
therefore a $G$--tree, with $n+1$ orbits of vertices and $n$ orbits of
edges.  Each edge orbit is free, one of the vertex orbits is free, and
there is one vertex orbit of type $G/G_i$ for each $1\leq i\leq n$.
Whenever a finite group acts on a
tree, it has a fixed point.  (To see this, take the finite subtree 
spanning an orbit, and peel off orbits of `leaves' until the remainder 
is fixed.)  Since the stabilizer of each edge
is trivial, it follows that each non-trivial finite subgroup of $G$
must fix exactly one vertex of the tree.  This implies that each
non-trivial finite subgroup of $G$ is conjugate to a subgroup of
exactly one of the $G_i$, and that two finite subgroups of $G_i$ are
conjugate in $G$ if and only if they were already conjugate in $G_i$.
\end{proof} 

\begin{proof}[Proof of Theorem~\ref{freeprod}] 
Let ${\cal Q}=\{Q_1,\ldots,Q_n\}$ be a finite list of isomorphism
types of finite groups not of prime power order.  For each $Q_i$, let
$G_i=H_i\semi Q_i$ be a group as in Theorem~\ref{ivf}.  Let $G$ be 
$G_1*\cdots * G_n$, the free product of the $G_i$.  By
Lemma~\ref{prodlem}, the group $G$ is of type $VF$, contains
infinitely many conjugacy classes of subgroup isomorphic to each
$Q_i$, and contains finitely many conjugacy classes of finite
subgroups of all other isomorphism types.
\end{proof} 

\begin{theorem} \label{fpq2} Let $Q$ be a non-trivial finite group.  
There exists a group $G=H\semi Q$ of type $FP$ over $\qq$, containing
infinitely many conjugacy classes of subgroups isomorphic to $Q$ and 
finitely many conjugacy classes of other finite subgroups.
Furthermore, $H$ is torsion-free, has rational cohomological
dimension at most~4 and has integral cohomological dimension at most~5.  
\end{theorem} 

\begin{proof} By Theorem~\ref{qacyc} there is a 3--dimensional
$\qq$--acyclic simplicial complex $L$ with a cocompact $\zz\times
Q$--action such that all stabilizers are finite, $L^Q=\emptyset$ and
$L^P\neq \emptyset$ if $P<Q$.  Take a flag triangulation of $L$, 
and consider the Bestvina--Brady group $H_L$.  By Theorem~\ref{type},
the semi-direct product $H=H_L\semi \zz$ is $FP$ over $\qq$.  By
Theorem~\ref{conjclass} the group $G=H_L\semi(\zz\times Q)$ contains
infinitely many conjugacy classes of subgroups isomorphic to $Q$ and 
finitely many conjugacy classes of other finite subgroups.  The
rational cohomological dimension of $H_L$ is at most the dimension of
the $\qq$--acyclic space $Y$ appearing in Section~\ref{bebr}, which is 
equal to the dimension of $L$, and the integral cohomological
dimension of $H_L$ is at most the dimension of the space $X_L$, which
is one more than the dimension of $L$.  The cohomological dimension of 
$H_L\semi \zz$ over any ring is at most one more than the
cohomological dimension of $H_L$ over the same ring. 
\end{proof} 

\begin{proof}[Proof of Theorem~\ref{freepq}] For each $Q_i\in Q$,
Theorem~\ref{fpq2} gives a group $G_i$ of type $FP$ over $\qq$ 
containing infinitely many conjugacy classes of subgroup isomorphic to
$Q$ and only finitely many conjugacy classes of other finite
subgroups.  By Lemma~\ref{prodlem}, the free product $G=G_1*\cdots
*G_n$ is $FP$ over $\qq$, contains infinitely many conjugacy classes
of subgroups isomorphic to each $Q_i\in \cal Q$, and contains finitely
many conjugacy classes of all other finite subgroups.  
\end{proof}

\begin{remark}  One difference between Theorems
\ref{freeprod}~and~\ref{freepq} is that each of the groups constructed
in Theorem~\ref{freepq} has virtual cohomological dimension at most~5,
whereas the virtual cohomological dimensions of the groups constructed
in Theorem~\ref{freeprod} seem to depend on the list $\cal Q$.  We do 
not know whether this necessarily happens, but the following proposition 
may be relevant.  
\end{remark} 

\begin{proposition} \label{dimprop}
Suppose that $G$ contains infinitely many conjugacy
classes of subgroup isomorphic to $SL_n(\ff_p)$, and that $G$ acts 
cocompactly with finite stabilizers on a mod-$p$--acyclic simplicial
complex $X$.  Then $X$ must have dimension at least $n-1$.  
\end{proposition} 

\begin{proof} There are only finitely many orbits in $X$, and hence
only finitely many conjugacy classes of subgroup of $G$ can fix some
point of $X$.  It follows that there is a subgroup isomorphic to
$SL_n(\ff_p)$ that has no fixed point, and we may apply
Theorem~\ref{slnp} to deduce the required result.  
\end{proof} 

\begin{remark} \label{qremark} 
If $G$ is virtually torsion-free and acts cocompactly with finite
stablizers on a contractible simplicial complex $X$, then $G$ is of
type $VF$.  It seems to be unknown whether every group of type 
$VF$ admits such an action.  It also seems to be unknown whether 
every group of type $FL$ over a prime field $F$ admits a free
cocompact action on an $F$--acyclic simplicial complex $X$.  
If $F$ is not assumed to be a prime field, then there are 
counterexamples.  In~\cite{epdn} we exhibited a group which is $FL$ over 
$\cc$ but which is not $FL$ over $\rr$.  This group cannot admit a 
cocompact free action on any 
$\cc$--acyclic simplicial complex $X$.  
\end{remark}

We conclude this section with a brief discussion of the Grothendieck
group $K_0(\qq G)$ of finitely generated projective modules for $\qq
G$ and its connection with conjugacy classes of elements of finite 
order in $G$.  First, we recall the definition of the Hattori--Stallings
trace~\cite{bass}.  

For any ring $R$, let $T(R)$ denote the quotient of $R$
by the additive subgroup generated by commutators of the form
$rs-sr$ for $r,s\in R$.  For a square matrix $A$ with coefficients
in $R$, the Hattori--Stallings trace $\tr(A)$ is the element of
$T(R)$ defined as the equivalence class containing the sum of the 
diagonal entries of $A$.  As an element of $T(R)$, this satisfies the 
usual trace condition $\tr(AB)=\tr(BA)$ for any matrices $A$ and $B$.  

Now suppose that $P$ is a finitely generated projective $R$ module, 
and that $P$ is isomorphic to a summand of $R^n$.  Pick an idempotent
$n\times n$ matrix $e_P$ whose image is isomorphic to $P$.  The 
Hattori--Stallings rank of $P$ is defined to be $\tr(e_P)$.  It may be
shown that this is independent of the choice of $n$ and $e_P$.  
The Hattori--Stallings rank defines a group homomorphism from $K_0(R)$
to $T(R)$.  

\begin{theorem}\label{bass}
For any group $G$, there is a subgroup of $K_0(\qq G)$ 
which is free abelian of rank equal to the number of conjugacy classes
of finite cyclic subgroups of $G$.  
\end{theorem} 

\begin{proof} For the group algebra $\qq G$, the group $T(\qq G)$ 
is the $\qq$--vector space with basis the conjugacy classes of elements
of $G$.  For any finite cyclic subgroup $C\leq G$, define an element 
$e_C\in \qq G$ by 
$$e_C= \frac{1}{|C|}\sum_{g\in C} g.$$ 
The element $e_C$ is an idempotent, and the $\qq G$--module $P_C$ defined by 
$P_C=\qq Ge_C$ is a projective $\qq G$--module.  With respect to the 
basis for $T(\qq G)$ given by the conjugacy classes of elements of
$G$, the non-zero coefficients in the Hattori--Stallings trace for 
$e_C$ are those corresponding to elements of $C$.  If $C_1,\ldots,
C_n$ are pairwise non-conjugate finite cyclic subgroups of $G$, it 
follows that the Hattori--Stallings traces $e_{C_1},\ldots,e_{C_n}$ are 
linearly independent.  It follows that the projectives of the form 
$P_C$ generate a subgroup of $K_0(\qq G)$ which is free abelian of 
rank equal to the number of conjugacy classes of finite cyclic
subgroups of $G$.  
\end{proof} 

\begin{corollary} There are groups $G$ of type $VF$ for which 
$K_0(\qq G)$ is not finitely generated.  
\end{corollary} 

\begin{proof} Apply Theorem~\ref{bass} to the groups with infinitely
  many conjugacy classes of finite cyclic subgroups constructed in 
Theorem~\ref{ivf}. 
\end{proof} 

\section{Other properties of the groups} 

Suppose that $Q$ is a group of automorphisms of a finite flag complex
$L$ with~$n$ vertices.  It is shown in~\cite{vfg} that in this case 
the group $G_L\semi Q$ is isomorphic to a subgroup of the special
linear group $SL_{2n}(\zz)$.  We do not know whether the groups 
$G_L\semi\Gamma$, for infinite $L$, are linear.  Residual finiteness 
however is easier to establish.  

\begin{lemma} \label{flagq}
Suppose that $\Gamma$ is residually finite and that 
$\Gamma$ acts cocompactly and with finite stabilizers on a flag
complex $L$.  There is a finite-index normal subgroup $\Gamma'$ such
that for any $\Gamma''\leq \Gamma'$, 
the quotient $L'=L/\Gamma''$ is a flag simplicial complex.  
\end{lemma} 

\begin{proof} There are finitely many conjugacy classes of simplex
stabilizer in $L$, and each simplex stabilizer is finite.  It follows
that there is a finite-index normal subgroup $\Gamma_1$ of $\Gamma$
that acts freely on $L$.  Since $L$ is locally finite, there are only
finitely many $\Gamma_1$--orbits of paths of length 1, 2 and 3 in the 
1--skeleton of $L$.  Hence we may pick $\Gamma_2$ of finite-index in 
$\Gamma_1$, such that no two points in the same $\Gamma_2$--orbit are 
joined by an edge path of length less than four.  We claim that we
may take $\Gamma'=\Gamma_2$.  

If $\Gamma''$ is any subgroup of $\Gamma_2$, then there is no edge
path of length less than four between any two vertices in the same 
$\Gamma''$--orbit.  In particular, there can be no loops in
$L/\Gamma'$.  Hence every simplex of $L$ maps injectively to a 
subspace of $L/\Gamma''$.  
There can be no double edges in $L/\Gamma''$, since that would give 
rise to an edge path of length two between vertices in the same 
$\Gamma''$--orbit.  Thus the 1--skeleton of $L/\Gamma''$ is a simplicial
complex.  

Now suppose that 
$\bar v_0,\ldots,\bar v_n$ are a mutually adjacent set of vertices 
of $L/\Gamma''$, and let $v_0$ be a lift of $\bar v_0$.  There exists 
a unique lift $v_i$ of each $\bar v_i$ that is adjacent to $v_0$.  
For each $i\neq j$, there exists a unique $g\in \Gamma_2$ so that 
$v_i$ is adjacent to $gv_j$.  But if $g\neq e$, then the path
$(v_j,v_0,v_i,gv_j)$ gives rise to a contradiction.  Thus the $v_i$ are all 
adjacent to each other, and so there is a simplex $\sigma$ 
of $L$ with vertex set $v_0,\ldots,v_n$.  
It follows that the quotient $L/\Gamma''$ contains a 
simplex $\bar \sigma$ spanning each complete subgraph of its
1--skeleton.  Suppose that $\bar\sigma'$ is any simplex of 
$L/\Gamma''$ spanning the same complete subgraph as $\bar \sigma$.  
There is a unique lift $\sigma'$ of $\bar \sigma'$ containing $v_0$.  
If $\sigma'\neq \sigma$, then there exists $i$ and $g\neq e$ so that 
$gv_i$ is a vertex of $\sigma'$.  But then there is an edge path 
of length 2 from $v_i$ to $gv_i$.  Hence any finite full subgraph of
the 1--skeleton of $L/\Gamma''$ is spanned by a unique simplex, and 
so $L/\Gamma''$ is a flag complex.  
\end{proof} 

\begin{theorem} \label{resfin} Let $\Gamma$ be residually finite and 
let $\Gamma$ act cocompactly and with finite stabilizers on a flag 
complex $L$.  Then the group $G_L\semi \Gamma$ is also residually finite.  
\end{theorem} 

\begin{proof} Let $g$ be a non-identity element of $G_L\semi\Gamma$.   
Since $(G_L\semi\Gamma)/G_L$ is isomorphic to $\Gamma$, it suffices to
consider the case when $g\in G_L$.  Let $K$ be a finite full
subcomplex of $L$ (ie, a subcomplex containing as many simplices as
possible given its 0--skeleton) such that $g$ is in the subgroup
generated by the vertices of $K$, and let $J$ be a finite full
subcomplex of $L$ containing $K$ and every vertex adjacent to a vertex
of $K$.  Let $\Gamma'$ be a finite-index subgroup of $\Gamma$ as in
Lemma~\ref{flagq}, and let $\Gamma''$ be a finite-index normal
subgroup of $\Gamma$ contained in $\Gamma'$ such that any two vertices
of $J$ lie in distinct $\Gamma''$--orbits.  Now $M=L/\Gamma''$ is a
finite flag complex, and $K$ maps to a full subcomplex of
$L/\Gamma''$.

The group $\Gamma/\Gamma''$ acts on $M$, and $g$ has non-trivial image
under the homomorphism $G_L\semi\Gamma\rightarrow
G_M\semi(\Gamma/\Gamma'')$.  Since this group is isomorphic to a
subgroup of $SL_{2n}(\zz)$, where $n$ is the number of vertices of
$M$ (see corollary~8 of \cite{vfg}), it follows that there is a finite
quotient of $G_L\semi\Gamma$ in which the image of $g$ is non-zero.  
\end{proof} 

In the special case when $\Gamma=\zz$ (which is the main case used
earlier in the paper), we shall show how to describe 
the group $G_L\semi\Gamma$ as the fundamental
group of a finite locally CAT(0) cube complex.  First we present two
lemmas concerning right-angled Artin groups. 

\begin{lemma}\label{artsub} Let $N$ be a full subcomplex of a flag
  complex $M$.  The inclusion $i\co N\rightarrow M$ induces a split 
injection $G_N\rightarrow G_M$.  
\end{lemma} 

\begin{proof} The quotient of $G_M$ by the subgroup generated by the 
vertices of $M-N$ is naturally isomorphic to $G_N$.  
\end{proof} 

\begin{lemma} \label{artgps}
Let the flag complex $K$ be expressed as $K= L\cup M$,
where $L$ and $M$ are full subcomplexes with $N=L\cap M$.  Then the 
group homomorphisms induced by the inclusion of each subcomplex in $K$ 
induce an isomorphism $G_L*_{G_N}G_M\rightarrow G_K$.  
\end{lemma} 

\begin{proof} 
Immediate from the presentations of the groups, given the result of 
Lemma~\ref{artsub}.  
\end{proof} 

\begin{theorem} \label{fpamalgam}
Let $\Gamma$ be an infinite cyclic group generated by
$\gamma$, let $\Gamma$ act on the flag complex $L$, and let $M$ be 
a `fundamental domain' for $\Gamma$ in the sense that $L=\bigcup_i
\gamma^i M$.  Define subcomplexes $N_0$ and $N_1$ by 
$$N_0= \gamma^{-1}M \cap M,\qquad N_1= M\cap \gamma M.$$ 
Then $G_L\semi \Gamma$ is isomorphic to the HNN--extension 
$G_M*_{G_{N_0}=G_{N_1}}$.  (In this HNN--extension, the base group is 
$G_M$, and the stable letter conjugates the subgroup $G_{N_0}$ to 
the subgroup $G_{N_1}$ by the map induced by 
$\gamma\co N_0\rightarrow N_1$.)  
\end{theorem}

\begin{proof} Let $t$ denote the stable letter in the HNN--extension,
and  consider the homomorphism $\phi$ from the HNN--extension to $\zz$ that 
sends $t$ to $1$ and sends each element of $G_M$ to $0$.  The kernel
of $\phi$ is an infinite free product with amalgamation: 
$$ \cdots *G_{-2}*_{H_{-1}}G_{-1}*_{H_0}G_0*_{H_1}
G_1*_{H_2}G_2*\cdots,$$ 
where $G_i$ denotes $t^{i}G_Mt^{-i}$, and $H_i$ denotes
$t^{i}G_{N_0}t^{-i}$.  If we define $M_i= \gamma^iM$ and
$N_i=\gamma^iN_0$, there is an isomorphism $\psi_i\co G_i\rightarrow
G_{M_i}$ defined as the composite
$$G_i \mapright{c(t^{-i})} G_0 \mapright{1} G_M \mapright{c(\gamma^i)} 
G_{M_i}$$ 
of conjugation by $t^{-i}$ followed by the identification of $G_0$ and 
$G_M$, followed by conjugation by $\gamma^i$.  Each of $\psi_i$ and 
$\psi_{i-1}$ induces an isomorphism from $H_i$ to $H_{N_i}$, and these
two are the same isomorphism.  The $\psi_i$ therefore fit together to
make an isomorphism from $\ker(\phi)$, described as an infinite free
product with amalgamation, to the following infinite free product 
with amalgamation:  
$$ \cdots *G_{M_{-2}}*_{H_{N_{-1}}}G_{M_{-1}}*_{H_{N_0}}
G_{M_0}*_{H_{N_1}}G_{M_1}*_{H_{N_2}}G_{M_2}*\cdots.$$ 
Furthermore, this isomorphism is equivariant for the $\zz$--actions
given by conjugation by powers of $t$ and $\gamma$.  
By Lemma~\ref{artgps}, the inclusions of the $G_{M_i}$ in $G_L$ induce
a $\Gamma$--equivariant isomorphism between the second free product 
with amalgamation and $G_L$.  Hence we obtain an isomorphism 
$G_M*_{G_{N_0}=G_{N_1}}\rightarrow G_L\semi \Gamma$ as required.  
\end{proof} 

\begin{corollary} \label{cubecclem}
Under the hypotheses of Theorem~\ref{fpamalgam}, 
the group $G_L\semi\Gamma$ is the fundamental group of a finite 
locally {\rm CAT(0)} cube complex.  
\end{corollary} 

\begin{proof} For any flag complex $K$, let $Y_K$ denote the explicit
model for the classifying space $BG_K$ described in
Section~\ref{bebr}, so that $Y_K=X_K/G_K$.  The naturality properties
of this construction are such that $Y_{N_0}$ and $Y_{N_1}$ are
subcomplexes of $Y_M$.  We construct a model $Z$ for $B(G_L\semi\Gamma)$ 
from $Y_M$ and $Y_{N_0}\times I$ by identifying 
$\{0\}\times Y_{N_0}$ with $Y_{N_0}\subseteq Y_M$ via the identity map
and identifying $\{1\}\times Y_{N_0}$ with $Y_{N_1}\subseteq Y_M$ via
the action of $\gamma$ which gives an isomorphism from $N_0$ to
$N_1$.  

The space $Z$ as above is a model for $B(G_L\semi \Gamma)$.  
To see that $Z$ has the structure of a locally CAT(0) cube complex,
one may either quote a gluing lemma (such as in \cite{BrHa},
proposition II.11.13), or
one may show that the link of the unique vertex in $Z$ is a flag
complex, which suffices by Gromov's lemma (\cite{BrHa}, theorem II.5.18).  
For any flag complex $K$, the link of the unique vertex in 
$Y_K$ is a flag complex $S(K)$, which is a sort of `double' of $K$:  
each vertex $v$ of $K$ corresponds to two vertices $v'$, $v''$ of 
$S(K)$, and a set of vertices of $S(K)$ is the vertex set of an 
$n$--simplex in $S(K)$ if and only if its image in the vertex set 
of $K$ is the vertex set of an $n$--simplex.  (For example, in the 
case when $K$ is a 2--simplex, $S(K)$ is the boundary of an
octahedron.)  The link of the vertex in $Z$ is isomorphic to 
$S(M)$ with a cone attached to each of the subspaces $S(N_0)$ and
$S(N_1)$, and hence it is a flag complex.  
\end{proof} 

\begin{corollary} Each of the groups $G_L\semi(\zz\times Q)$
constructed in Section~\ref{examps} acts cocompactly with finite 
stabilizers on some {\rm CAT(0)} cube complex.  In particular, there is a
model for the universal space for proper actions of 
$G_L\semi(\zz\times Q)$ which has finitely many orbits of cells.  
\end{corollary} 

\begin{proof} Take a finite `fundamental domain', $M'$, for the action
of $\zz$ on $L$ (as in the statement of Theorem~\ref{fpamalgam}).  In 
case $M'$ is not $Q$--invariant, replace $M'$ by $M= \bigcup_{q\in Q}
qM'$.  For this choice of $M$, there is a base-point preserving 
cellular $Q$--action on $Z$, the model for $B(G_L\semi \zz)$ constructed in 
Corollary~\ref{cubecclem}.  This induces the required action of
$G_L\semi(\zz\times Q)$ on the universal cover of $Z$.  Whenever a 
group $H$ acts with finite stabilizers on a CAT(0) cube complex, that 
space is a model for the universal space for proper actions of $H$
\cite{vfg}.  
\end{proof}

\end{document}